# Common Factors in Fraction-Free Matrix Decompositions

Johannes Middeke, David J. Jeffrey and Christoph Koutschan

**Abstract.** We consider $LU$ and $QR$ matrix decompositions using exact computations. We show that fraction-free Gauß–Bareiss reduction leads to triangular matrices having a non-trivial number of common row factors. We identify two types of common factors: systematic and statistical. Systematic factors depend on the reduction process, independent of the data, while statistical factors depend on the specific data. We relate the existence of row factors in the $LU$ decomposition to factors appearing in the Smith–Jacobson normal form of the matrix. For statistical factors, we identify some of the mechanisms that create them and give estimates of the frequency of their occurrence. Similar observations apply to the common factors in a fraction-free $QR$ decomposition. Our conclusions are tested experimentally.

## 1. Introduction

Although known earlier, the fraction-free method for exact matrix computations became well known because of its application by Bareiss [1] to the solution of a system on integer equations $AX = B$. He implemented fraction-free Gaussian elimination of the augmented matrix $[A\ B]$, and retained integer computations until a final division step. Since, in linear algebra, equation solving is related to the matrix factorizations $LU$ and $QR$, it was natural that fraction-free methods would be extended later to those topics. The extensions required that the forms of the factorizations be modified from their floating-point counterparts in order to retain purely integer data. The first proposals were based on inflating the initial data until all divisions were guaranteed exact, see for example Lee and Saunders [12]; Nakos et al. [14]; Corless and Jeffrey [4]. This strategy, however, led to the entries in the $L$ and $U$ matrices becoming very large, and an alternative form was presented in Zhou and Jeffrey [18]. Fraction-free Gram–Schmidt orthogonalization and $QR$ factorization were similarly studied by Erlingsson et al. [5]; Zhou and Jeffrey [18]. Further extensions have addressed fraction-free full-rank factoring of non-invertible matrices and fraction-free computation of the Moore-Penrose inverse [10]. More generally, applications exist in areas such as the Euclidean algorithm, and the Berlekamp–Massey algorithm [11].

More general domains are possible, and here we consider matrices over a principal ideal domain $\mathbb{D}$. For the purpose of giving illustrative examples and conducting computational experiments, matrices over $\mathbb{Z}$ and $\mathbb{Q}[x]$ are used, because these domains are well established and familiar to readers. We emphasize, however, that the methods here apply for all principal ideal domains, as opposed to methods that target specific domains, such as Giesbrecht and Storjohann [7]; Pauderis and Storjohann [16].

The shift from equation solving to matrix factorization has the effect of making visible the intermediate results that are deleted in the original Bareiss implementation. Because of this, it becomes

---
J. M. was supported in part by the Austrian Science Fund (FWF): SFB50 (F5009-N15).
C. K. was supported by the Austrian Science Fund (FWF): P29467-N32 and F5011-N15.

apparent that the columns and rows of the $L$ and $U$ matrices frequently contain common factors, which otherwise pass unnoticed. We consider here how these factors arise, and what consequences there are for the computations.

Our starting point is a fraction-free form for $LU$ decomposition [10]: given a matrix $A$ over $\mathbb{D}$,

$$A = P_r L D^{-1} U P_c,$$

where $L$ and $U$ are lower and upper triangular matrices, respectively, $D$ is a diagonal matrix, and the entries of $L$, $D$, and $U$ are from $\mathbb{D}$. The permutation matrices $P_r$ and $P_c$ ensure that the decomposition is always a full-rank decomposition, even if $A$ is rectangular or rank deficient; see section 2. The decomposition is computed by a variant of Bareiss's algorithm [1]. In section 6, the $LD^{-1}U$ decomposition also is the basis of a fraction-free $QR$ decomposition.

The key feature of Bareiss's algorithm is that it creates factors which are common to every element in a row, but which can then be removed by exact divisions. We refer to such factors, which appear predictably owing to the decomposition algorithm, as "systematic factors". There are, however, other common factors which occur with computable probability, but which depend upon the particular data present in the input matrix. We call such factors "statistical factors". In this paper we discuss the origins of both kinds of common factors and show that we can predict a nontrivial proportion of them from simple considerations.

Once the existence of common factors is recognized, it is natural to consider what consequences, if any, there are for the computation, or application, of the factorizations. Some consequences we shall consider include a lack of uniqueness in the definition of the $LU$ factorization, and whether the common factors add significantly to the sizes of the elements in the constituent factors. This in turn leads to questions regarding the benefits of removing common factors, and what computational cost is associated with such benefits.

A synopsis of the paper is as follows. After recalling Bareiss's algorithm, the $LD^{-1}U$ decomposition, and the algorithm from Jeffrey [10] in section 2, we establish, in section 3, a relation between the systematic common row factors of $U$ and the entries in the Smith–Jacobson normal form of the same input matrix $A$. In section 4 we propose an efficient way of identifying some of the systematic common row factors introduced by Bareiss's algorithm; these factors can then be easily removed by exact division. In section 5 we present a detailed analysis concerning the expected number of statistical common factors in the special case $\mathbb{D} = \mathbb{Z}$, and we find perfect agreement with our experimental results. We conclude that the factors make a measurable contribution to the element size, but they do not impose a serious burden on calculations.

In section 6 we investigate the $QR$ factorization. In this context, the orthonormal $Q$ matrix used in floating point calculations is replaced by a $\Theta$ matrix, which is left-orthogonal, i.e. $\Theta^t \Theta$ is diagonal, but $\Theta \Theta^t$ is not. We show that, for a square matrix $A$, the last column of $\Theta$, as calculated by existing algorithms, is subject to an exact division by the determinant of $A$, with a possibly significant reduction in size.

Throughout the paper, we employ the following notation. We assume, unless otherwise stated, that the ring $\mathbb{D}$ is an arbitrary principal ideal domain. We denote the set of all $m$-by-$n$ matrices over $\mathbb{D}$ by $\mathbb{D}^{m \times n}$. We write $\mathbf{1}_n$ for the $n$-by-$n$ identity matrix and $\mathbf{0}_{m \times n}$ for the $m$-by-$n$ zero matrix. We shall usually omit the subscripts if no confusion is possible. For $A \in \mathbb{D}^{m \times n}$ and $1 \leq i \leq m$, $A_{i,*}$ is the $i^{\text{th}}$ row of $A$. Similarly, $A_{*,j}$ is the $j^{\text{th}}$ column of $A$ for $1 \leq j \leq n$. If $1 \leq i_1 < i_2 \leq m$ and $1 \leq j_1 < j_2 \leq n$, we use $A_{i_1 \ldots i_2, j_1 \ldots j_2}$ to refer to the submatrix of $A$ made up from the entries of the rows $i_1$ to $i_2$ and the columns $j_1$ to $j_2$. Given elements $a_1, \ldots, a_n \in \mathbb{D}$, with $\operatorname{diag}(a_1, \ldots, a_n)$ we refer to the diagonal matrix that has $a_j$ as the entry at position $(j, j)$ for $1 \leq j \leq n$. We will use the same notation for block diagonal matrices.

We denote the set of all column vectors of length $m$ with entries in $\mathbb{D}$ by $\mathbb{D}^m$ and that of all row vectors of length $n$ by $\mathbb{D}^{1 \times n}$. If $\mathbb{D}$ is a unique factorization domain and $v = (v_1, \ldots, v_n) \in \mathbb{D}^{1 \times n}$, then we set $\gcd(v) = \gcd(v_1, \ldots, v_n)$. Moreover, with $d \in \mathbb{D}$ we write $d \mid v$ if $d \mid v_1 \wedge \ldots \wedge d \mid v_n$ (or, equivalently, if $d \mid \gcd(v)$). We also use the same notation for column vectors.

We will sometimes write column vectors $w \in \mathbb{D}^m$ with an underline $\underline{w}$ and row vectors $v \in \mathbb{D}^{1 \times n}$ with an overline $\overline{v}$ if we want to emphasize the specific type of vector.

## 2. Bareiss's Algorithm and the $LD^{-1}U$ Decomposition

For the convenience of the reader, we start by recalling Bareiss's algorithm from Bareiss [1] as well as the $LD^{-1}U$ decomposition from Jeffrey [10].

Let $\mathbb{D}$ be an integral domain[1], and let $A \in \mathbb{D}^{m \times n}$ be a matrix and $b \in \mathbb{D}^m$ be a vector. We are interested in computing the exact solution(s) for the system $Ax = b$ over the quotient field of $\mathbb{D}$. Following standard methods from linear algebra, we compute a (variation of the) $LU$ decomposition for $A$. Suppose that

$$A = \begin{pmatrix} u & \cdots \\ v & \cdots \\ \vdots & \ddots \end{pmatrix},$$

where we assume that $u$ and $v$ are not zero. Choosing $u$ as our pivot in the usual Gaussian elimination, we would subtract $v/u$ times the first row from the second as seen in, for example, Geddes et al. [6, Section 9.2]. The division, however, forces us to work in the quotient field of $\mathbb{D}$, which is computationally expensive. One way to avoid this is to use cross-multiplication (see Geddes et al. [6, Section 9.2] where it is called division-free elimination): In order to eliminate $v$ from the second row, first multiply the second row by $u$ and then subtract $v$ times the first row. While this avoids dealing with fractions, when repeating this process for the entire elimination the matrix entries grow in size exponentially with the number of rows, which will also need more computation time (and take up more memory). If $\mathbb{D}$ is a unique factorization domain, it is possible to keep the sizes of the entries as small as possible by dividing each row by its greatest common divisor after the cross-multiplication. The drawback is that computing these greatest common divisors in all steps is again computationally expensive.

Bareiss [1] observed that after carrying out cross-multiplication for two iterations (that is, eliminating the entries from the first column except the first and eliminating the entries from the second column except the first two) every entry in the third row and below is divisible by $u$. More generally, if we denote the pivots during the elimination by $p_1, \ldots, p_r$ (where $r$ is the rank of $A$), then after $k \geq 2$ iterations every entry below the $k^{\text{th}}$ row is divisible by $p_{k-1}$. (In actual implementations one often additionally defines $p_{-1} = 1$ in order to prevent $k = 1$ from being a special case.) Bareiss's observation gives us the possibility of removing the factor $p_{k-1}$ after the $k^{\text{th}}$ iteration from the lower rows of the matrix by using (relatively) cheap exact division while at the same time avoiding the cost of computing the greatest common divisors. This makes Bareiss's algorithm a valuable choice for exact matrix computations (see, for example, Geddes et al. [6, Section 9.3] where it is called the single-step fraction-free elimination scheme). We will show below that we can extract even more factors at a moderate additional cost.

In Jeffrey [10] the idea of Bareiss's algorithm was used in order to obtain a fraction-free variant of the $LU$ factorization. We repeat the main result from that paper here.

**Theorem 1 (Jeffrey [10, Thm. 2]).** *A rectangular matrix $A$ with elements from an integral domain $\mathbb{D}$, having dimensions $m \times n$ and rank $r$, may be factored into matrices containing only elements from $\mathbb{D}$ in the form*

$$A = P_r L D^{-1} U P_c = P_r \begin{pmatrix} \mathcal{L} \\ \mathcal{M} \end{pmatrix} D^{-1} \begin{pmatrix} \mathcal{U} & \mathcal{V} \end{pmatrix} P_c$$

---

[1]Note that in this section we do not require $\mathbb{D}$ to be a principal ideal domain, but it suffices to assume that $\mathbb{D}$ is an integral domain.

where the permutation matrix $P_r$ is $m \times m$; the permutation matrix $P_c$ is $n \times n$; $\mathcal{L}$ is $r \times r$, lower triangular and has full rank:

$$\mathcal{L} = \begin{pmatrix} p_1 & 0 & \cdots & 0 \\ \ell_{21} & p_2 & \ddots & \vdots \\ \vdots & \vdots & \ddots & 0 \\ \ell_{r1} & \ell_{r2} & \cdots & p_r \end{pmatrix}$$

where the $p_i \neq 0$ are the pivots in a Gaussian elimination; $\mathcal{M}$ is $(m-r) \times r$ and could be null; $D$ is $r \times r$ and diagonal:

$$D = \operatorname{diag}(p_1, p_1 p_2, p_2 p_3, \ldots, p_{r-2} p_{r-1}, p_{r-1} p_r);$$

$\mathcal{U}$ is $r \times r$ and upper triangular, while $\mathcal{V}$ is $r \times (n-r)$ and could be null:

$$\mathcal{U} = \begin{pmatrix} p_1 & u_{12} & \cdots & u_{1r} \\ 0 & p_2 & \cdots & u_{2r} \\ \vdots & \ddots & \ddots & \vdots \\ 0 & \cdots & 0 & p_r \end{pmatrix}.$$

Inspecting the proof given in Jeffrey [10], it is possible to extract an algorithm for the computation of the $LD^{-1}U$ decomposition. Note that the algorithm can be seen as a variant of Bareiss's algorithm [1] which will yield the same $U$. The difference is that Jeffrey [10] also explains how to obtain $L$ and $D$ in a fraction-free way.

**Algorithm 2.** *($LD^{-1}U$ decomposition)*

**Input:.** *A matrix $A \in \mathbb{D}^{m \times n}$.*
**Output:.** *The $LD^{-1}U$ decomposition of $A$ as in Theorem 1.*

1. *Initialize $p_{-1} = 1$, $P_r = \mathbf{1}_m$, $L = \mathbf{0}_{m \times m}$, $U = A$ and $P_c = \mathbf{1}_n$.*
2. *For each $k = 1, \ldots, \min\{m, n\}$:*
   (a) *Find a non-zero pivot $p_k$ in $U_{k\ldots m, k\ldots n}$ and bring it to position $(k,k)$ recording the row and column swaps in $P_r$ and $P_c$. Also apply the row swaps to $L$ accordingly. If no pivot is found, then set $r = k$ and exit the loop.*
   (b) *Set $L_{k,k} = p_k$ and $L_{i,k} = U_{i,k}$ for $i = k+1, \ldots, m$.*
   *Then eliminate the entries in the $k^{\text{th}}$ column and below the $k^{\text{th}}$ row in $U$ by cross-multiplication; that is, for $i > k$ set $U_{i,*}$ to $p_k U_{i,*} - U_{ik} U_{k,*}$.*
   (c) *Perform division by $p_{k-1}$ on the rows beneath the $k^{\text{th}}$ in $U$; that is, for $i > k$ set $U_{i,*}$ to $U_{i,*}/p_{k-1}$. Note that the divisions will be exact.*
3. *If $r$ is not set yet, set $r = \min\{m, n\}$.*
4. *If $r < m$, then trim the last $m - r$ columns from $L$ as well as the last $m - r$ rows from $U$.*
5. *Set $D = \operatorname{diag}(p_1, p_1 p_2, \ldots, p_{r-1} p_r)$.*
6. *Return $P_r$, $L$, $D$, $U$, and $P_c$.*

The algorithm does not specify the choice of pivot in step 2a. Conventional wisdom (see, for example, Geddes et al. [6]) is that in exact algorithms choosing the smallest possible pivot (measured in a way suitable for $\mathbb{D}$) will lead to the smallest output sizes. We have been able to confirm this experimentally in Middeke and Jeffrey [13] for $\mathbb{D} = \mathbb{Z}$ where size was measured as the absolute value. In step 2c the divisions are guaranteed to be exact. Thus, an implementation can use more efficient procedures for this step if available (for example, for big integers using `mpz_divexact` in the GMP library which is based on Jebelean [9] instead of regular division).

One of the goals of the present paper is to discuss improvements to the decomposition explained above. Throughout this paper we shall use the term $LD^{-1}U$ *decomposition* to mean exactly the decomposition from Theorem 1 as computed by Algorithm 2. For the variations of this decomposition we introduce the following term:

**Definition 3 (Fraction-Free $LU$ Decomposition).** For a matrix $A \in \mathbb{D}^{m \times n}$ of rank $r$ we say that $A = P_r L D^{-1} U P_c$ is a *fraction-free LU decomposition* if $P_r \in \mathbb{D}^{m \times m}$ and $P_c \in \mathbb{D}^{n \times n}$ are permutation matrices, $L \in \mathbb{D}^{m \times r}$ has $L_{ij} = 0$ for $j > i$ and $L_{ii} \neq 0$ for all $i$, $U \in \mathbb{D}^{r \times n}$ has $U_{ij} = 0$ for $i > j$ and $U_{ii} \neq 0$ for all $i$, and $D \in \mathbb{D}^{r \times r}$ is a diagonal matrix (with full rank).

We will usually refer to matrices $L \in \mathbb{D}^{m \times r}$ with $L_{ij} = 0$ for $j > i$ and $L_{ii} \neq 0$ for all $i$ as *lower triangular* and to matrices $U \in \mathbb{D}^{r \times n}$ with $U_{ij} = 0$ for $i > j$ and $U_{ii} \neq 0$ for all $i$ as *upper triangular* even if they are not square.

As mentioned in the introduction, Algorithm 2 does result in common factors in the rows of the output $U$ and the columns of $L$. In the following sections, we will explore methods to explain and predict those factors. The next result asserts that we can cancel all common factors which we find from the final output. This yields a fraction-free $LU$ decomposition of $A$ where the size of the entries of $U$ (and $L$) are smaller than in the $LD^{-1}U$ decomposition.

**Corollary 4.** *Given a matrix $A \in \mathbb{D}^{m \times n}$ with rank $r$ and its standard $LD^{-1}U$ decomposition $A = P_c L D^{-1} U P_c$, if $D_U = \mathrm{diag}(d_1, \ldots, d_r)$ is a diagonal matrix with $d_k \mid \gcd(U_{k,*})$, then setting $\hat{U} = D_U^{-1} U$ and $\hat{D} = D D_U^{-1}$ where both matrices are fraction-free we have the decomposition $A = P_c L \hat{D}^{-1} \hat{U} P_c$.*

*Proof.* By Theorem 1, the diagonal entries of $U$ are the pivots chosen during the decomposition and they also divide the diagonal entries of $D$. Thus, any common divisor of $U_{k,*}$ will also divide $D_{kk}$ and therefore both $\hat{U}$ and $\hat{D}$ are fraction-free. We can easily check that $A = P_c L D^{-1} D_U D_U^{-1} U = P_c L \hat{D}^{-1} \hat{U} P_c$. □

**Remark 5.** If we predict common column factors of $L$ we can cancel them in the same way. However, if we have already canceled factors from $U$, then there is no guarantee that $d \mid L_{*,k}$ implies $d \mid \hat{D}_{kk}$. Thus, in general we can only cancel $\gcd(d, \hat{D}_{kk})$ from $L_{*,k}$. The same holds *mutatis mutandis* if we cancel the factors from $L$ first.

It will be an interesting discussion for future research whether it is better to cancel as many factors as possible from $U$ or to cancel them from $L$.

## 3. LU and the Smith–Jacobson Normal Form

This section explains a connection between "systematic factors" (that is, common factors which appear in the decomposition due to the algorithm being used) and the Smith–Jacobson normal form. Given a matrix $A$ over a principal ideal domain $\mathbb{D}$, we study the decomposition $A = P_r L D^{-1} U P_c$. For simplicity, from now on we consider the decomposition in the form $P_r^{-1} A P_c^{-1} = L D^{-1} U$. The following theorem connecting the $LD^{-1}U$ decomposition with the Smith–Jacobson normal form can essentially be found in [1].

**Theorem 6.** *Let the matrix $A \in \mathbb{D}^{n \times n}$ have the Smith–Jacobson normal form $S = \mathrm{diag}(d_1, \ldots, d_n)$ where $d_1, \ldots, d_n \in \mathbb{D}$. Moreover, let $A = L D^{-1} U$ be an $L D^{-1} U$ decomposition of $A$ without permutations. Then for $k = 1, \ldots, n$*

$$d_k^* = \prod_{j=1}^{k} d_j \mid U_{k,*} \quad \text{and} \quad d_k^* \mid L_{*,k}.$$

**Remark 7.** The values $d_1^*, \ldots, d_n^*$ are known in the literature as the *determinantal divisors* of $A$.

*Proof.* According to [15, II.15], the diagonal entries of the Smith–Jacobson normal form are quotients of the determinantal divisors, i.e., $d_1^* = d_1$ and $d_k = d_k^*/d_{k-1}^*$ for $k = 2, \ldots, n$. Moreover, $d_k^*$ is the greatest common divisor of all $k$-by-$k$ minors of $A$ for each $k = 1, \ldots, n$. Thus, we only have to prove that the entries of the $k^{\text{th}}$ row of $U$ are $k$-by-$k$ minors of $A$. However, this follows from [6, Eqns (9.8), (9.12)], since the $k^{\text{th}}$ row of $U$ consists of the elements

$$U_{kj} = \det \begin{pmatrix} A_{1,1} & \cdots & A_{1,k-1} & A_{1,j} \\ \vdots & & \vdots & \vdots \\ A_{k,1} & \cdots & A_{k,k-1} & A_{k,j} \end{pmatrix} \quad \text{where} \quad j = 1, \ldots, n.$$

Similarly, following the algorithm in Jeffrey [10], we see that the columns of $L$ are just made up by copying entries from the columns of $U$ during the reduction. More precisely, the $k^{\text{th}}$ column of $L$ will have the entries $a_{1k}^{(k-1)}, \ldots, a_{nk}^{(k-1)}$, using the notation of Geddes et al. [6], and these are again just $k$-by-$k$ minors of $A$. Explicitly, the $j^{\text{th}}$ column of $L$ consists of the elements

$$L_{kj} = \det \begin{pmatrix} A_{1,1} & \cdots & A_{1,j} \\ \vdots & & \vdots \\ A_{j-1,1} & \cdots & A_{j-1,j} \\ A_{k,1} & \cdots & A_{k,j} \end{pmatrix}$$

where $k = 1, \ldots, n$. □

¿From Theorem 6, we obtain the following result.

**Corollary 8.** *The $k^{\text{th}}$ determinantal divisor $d_k^*$ can be removed from the $k^{\text{th}}$ row of $U$ (since it divides $D_{k,k}$ by Corollary 4) and also $d_{k-1}^*$ can be removed from the $k^{\text{th}}$ row of $L$ because $d_{k-1}^* \mid d_k^*$ and $d_j^*$ divides the $j^{\text{th}}$ pivot for $j = k-1, k$. Thus, $d_{k-1}^* d_k^* \mid D_{k,k}$.*

We give an example using the domain $\mathbb{Q}[x]$. Let $A$ be the polynomial matrix

$$\begin{pmatrix} -\frac{3}{2} & -x^3 + 5x^2 + 3x - \frac{9}{2} & x^2 + x & \frac{1}{2}x^3 - x^2 \\ -3 & -2x^3 + 10x^2 + 5x - 9 & 2x^2 + 2x & x^3 - 2x^2 \\ \frac{1}{2} & x^3 + \frac{3}{2} & 0 & -\frac{1}{2}x^3 \\ -\frac{1}{2} & -x - \frac{3}{2} & 0 & \frac{1}{2}x \end{pmatrix}.$$

The Smith–Jacobson normal form $S$ of $A$ is

$$\operatorname{diag}(1, x, x(x+1), x(x+1)(x-1))$$

and thus its determinantal divisors are $d_1^* = 1$, $d_2^* = x$, $d_3^* = x^2(x+1)$ and $d_4^* = x^3(x+1)^2(x-1)$. Computing the $LD^{-1}U$ decomposition of $A$ yields $A = LD^{-1}U$ where $L$ is

$$\begin{pmatrix} -\frac{3}{2} & 0 & 0 & 0 \\ -3 & \frac{3}{2}x & 0 & 0 \\ \frac{1}{2} & -x^3 - \frac{5}{2}x^2 - \frac{3}{2}x & \frac{1}{2}x^3 + \frac{1}{2}x^2 & 0 \\ -\frac{1}{2} & -\frac{1}{2}x^3 + \frac{5}{2}x^2 + 3x & -\frac{1}{2}x^3 - \frac{1}{2}x^2 & -\frac{1}{4}x^6 - \frac{1}{4}x^5 + \frac{1}{4}x^4 + \frac{1}{4}x^3 \end{pmatrix},$$

$D = \operatorname{diag}(-3/2, -9/4x, 3/4x^4 + 3/4x^3, -1/8x^9 - 1/4x^8 + 1/4x^6 + 1/8x^5)$, and $U$ is

$$\begin{pmatrix} -\frac{3}{2} & -x^3 + 5x^2 + 3x - \frac{9}{2} & x^2 + x & \frac{1}{2}x^3 - x^2 \\ 0 & \frac{3}{2}x & 0 & 0 \\ 0 & 0 & \frac{1}{2}x^3 + \frac{1}{2}x^2 & -\frac{1}{2}x^4 - \frac{1}{2}x^3 \\ 0 & 0 & 0 & -\frac{1}{4}x^6 - \frac{1}{4}x^5 + \frac{1}{4}x^4 + \frac{1}{4}x^3 \end{pmatrix}.$$

Computing the column factors of $L$ and the row factors of $U$ yields the list $1$, $x$, $x^2(x+1)$ and $x^3(x-1)(x+1)^2$, i.e., exactly the determinantal divisors. In general, there could be other factors as well.

## 4. Efficient Detection of Factors

When considering the output of Algorithm 2, we find an interesting relation between the entries of $L$ and $U$ which can be exploited in order to find "systematic" common factors in the $LD^{-1}U$ decomposition. Theorem 9 below predicts a divisor of the common factor in the $k^{\text{th}}$ row of $U$, by looking at just three entries of $L$. Likewise, we obtain a divisor of the common factor of the $k^{\text{th}}$ column of $L$ from three entries of $U$. As in the previous section, let $\mathbb{D}$ be a principal ideal domain.

**Theorem 9.** *Let $A \in \mathbb{D}^{m \times n}$ and let $P_r L D^{-1} U P_c$ be the $LD^{-1}U$ decomposition of $A$. Then*

$$\frac{\gcd(L_{k-1,k-1}, L_{k,k-1})}{\gcd(L_{k-1,k-1}, L_{k,k-1}, L_{k-2,k-2})} \;\Big|\; U_{k,*}$$

*and*

$$\frac{\gcd(U_{k-1,k-1}, U_{k-1,k})}{\gcd(U_{k-1,k-1}, U_{k-1,k}, U_{k-2,k-2})} \;\Big|\; L_{*,k}$$

*for $k = 2, \ldots, m-1$ (where we use $L_{0,0} = U_{0,0} = 1$ for $k=2$).*

*Proof.* Suppose that during Bareiss's algorithm after $k-1$ iterations we have reached the following state

$$A^{(k-1)} = \begin{pmatrix} T & \underline{*} & \underline{*} & * \\ \overline{0} & p & * & \overline{*} \\ \overline{0} & 0 & a & \overline{v} \\ \overline{0} & 0 & b & \overline{w} \\ \mathbf{0} & \underline{0} & \underline{*} & * \end{pmatrix},$$

where $T$ is an upper triangular matrix, $p, a, b \in \mathbb{D}$, $\overline{v}, \overline{w} \in \mathbb{D}^{1 \times n-k-1}$ and the other overlined quantities are row vectors and the underlined quantities are column vectors. Assume that $a \neq 0$ and that we choose it as a pivot. Continuing the computations we now eliminate $b$ (and the entries below) by cross-multiplication

$$A^{(k-1)} \leadsto \begin{pmatrix} T & \underline{*} & \underline{*} & * \\ \overline{0} & p & * & \overline{*} \\ \overline{0} & 0 & a & \overline{v} \\ \overline{0} & 0 & 0 & a\overline{w} - b\overline{v} \\ \mathbf{0} & \underline{0} & \underline{0} & * \end{pmatrix}.$$

Here, we can see that any common factor of $a$ and $b$ will be a factor of every entry in that row, i.e., $\gcd(a,b) \mid a\overline{w} - b\overline{v}$. However, we still have to carry out the exact division step. This leads to

$$A^{(k-1)} \leadsto \begin{pmatrix} T & \underline{*} & \underline{*} & * \\ \overline{0} & p & * & \overline{*} \\ \overline{0} & 0 & a & \overline{v} \\ \overline{0} & 0 & 0 & \frac{1}{p}(a\overline{w} - b\overline{v}) \\ \mathbf{0} & \underline{0} & \underline{0} & * \end{pmatrix} = A^{(k)}.$$

The division by $p$ is exact. Some of the factors in $p$ might be factors of $a$ or $b$ while others are hidden in $\overline{v}$ or $\overline{w}$. However, every common factor of $a$ and $b$ which is not also a factor of $p$ will still be a common factor of the resulting row. In other words,

$$\frac{\gcd(a,b)}{\gcd(a,b,p)} \;\Big|\; \frac{1}{p}(a\overline{w} - b\overline{v}).$$

In fact, the factors do not need to be tracked during the $LD^{-1}U$ reduction but can be computed afterwards: All the necessary entries $a$, $b$ and $p$ of $A^{(k-1)}$ will end up as entries of $L$. More precisely, we shall have $p = L_{k-2,k-2}$, $a = L_{k-1,k-1}$ and $b = L_{k,k-1}$.

Similar reasoning can be used to predict common factors in the columns of $L$. Here, we have to take into account that the columns of $L$ are made up from entries in $U$ during each iteration of the computation. □

As a typical example consider the matrix

$$A = \begin{pmatrix} 8 & 49 & 45 & -77 & 66 \\ -10 & -77 & -19 & -52 & 48 \\ 51 & 18 & -81 & 31 & 69 \\ -97 & -58 & 37 & 41 & 22 \\ -60 & 0 & -25 & -18 & -92 \end{pmatrix}.$$

This matrix has a $LD^{-1}U$ decomposition with

$$L = \begin{pmatrix} 8 & 0 & 0 & 0 & 0 \\ -10 & -126 & 0 & 0 & 0 \\ 51 & -2355 & 134076 & 0 & 0 \\ -97 & 4289 & -233176 & -28490930 & 0 \\ -60 & 2940 & -148890 & -53377713 & 11988124645 \end{pmatrix}$$

and with

$$U = \begin{pmatrix} 8 & 49 & 45 & -77 & 66 \\ 0 & -126 & 298 & -1186 & 1044 \\ 0 & 0 & 134076 & -414885 & 351648 \\ 0 & 0 & 0 & -28490930 & 55072620 \\ 0 & 0 & 0 & 0 & 11988124645 \end{pmatrix}.$$

Note that in this example pivoting is not needed, that is, we have $P_r = P_c = \mathbf{1}$. The method outlined in Theorem 9 correctly predicts the common factor 2 in the second row, the factor 3 in the third row and the factor 2 in the fourth row. However, it does not detect the additional factor 5 in the fourth row.

The example also provides an illustration of the proof of Theorem 6: The entry $-414885$ of $U$ at position $(3,4)$ is given by the determinant of the submatrix

$$\begin{pmatrix} 8 & 49 & -77 \\ -10 & -77 & -52 \\ 51 & 18 & 31 \end{pmatrix}$$

consisting of the first three rows and columns 1, 2 and 4 of $A$. In this particular example, however, the Smith–Jacobson Normal Form of the matrix $A$ is $\mathrm{diag}(1,1,1,1,11988124645)$ which does not yield any information about the common factors.

Given Theorem 9, one can ask how good this prediction actually is. Concentrating on the case of integer matrices, the following Theorem 10 shows that with this prediction we do find a common factor in roughly a quarter of all rows. Experimental data suggest a similar behavior for matrices containing polynomials in $\mathbb{F}_p[x]$ where $p$ is prime. Moreover, these experiments also showed that the prediction was able to account for 40.17% of all the common prime factors (counted with multiplicity) in the rows of $U$.[2]

**Theorem 10.** *For random integers $a, b, p \in \mathbb{Z}$ the probability that the formula in Theorem 9 predicts a non-trivial common factor is*

$$\mathrm{P}\Big(\frac{\gcd(a,b)}{\gcd(p,a,b)} \neq 1\Big) = 6\frac{\zeta(3)}{\pi^2} \approx 26.92\%.$$

*Proof.* The following calculation is due to Hare [8]; Winterhof [17]: First note that the probability that $\gcd(a,b) = n$ is $1/n^2$ times the probability that $\gcd(a,b) = 1$. Summing up all of these probabilities gives

$$\sum_{n=1}^{\infty} \mathrm{P}\big(\gcd(a,b) = n\big) = \sum_{n=1}^{\infty} \frac{1}{n^2} \mathrm{P}\big(\gcd(a,b) = 1\big) = \mathrm{P}\big(\gcd(a,b) = 1\big)\frac{\pi^2}{6}.$$

---

[2]This experiment was carried out with random square matrices $A$ of sizes between 5-by-5 and 125-by-125. We decomposed $A$ into $P_r L D^{-1} U P_c$ and then computed the number of predicted prime factors in $U$ and related that to the number of actual prime factors. We did not consider the last row of $U$ since this contains only the determinant.

As this sum must be 1, this gives that the $\mathrm{P}\bigl(\gcd(a,b)=1\bigr)=6/\pi^2$, and the $\mathrm{P}\bigl(\gcd(a,b)=n\bigr)=6/(\pi^2 n^2)$. Given that $\gcd(a,b)=n$, the probability that $n\mid c$ is $1/n$. So the probability that $\gcd(a,b)=n$ and that $\gcd(p,a,b)=n$ is $6/(\pi^2 n^3)$. So $\mathrm{P}\bigl(\gcd(a,b)/\gcd(p,a,b)=1\bigr)$ is

$$\sum_{n=1}^{\infty} \mathrm{P}\bigl(\gcd(a,b)=n \text{ and } \gcd(p,a,b)=n\bigr) = \sum_{n=1}^{\infty} \frac{6}{\pi^2 n^3} = 6\frac{\zeta(3)}{\pi^2}. \qquad \square$$

There is another way in which common factors in integer matrices can arise. Let $d$ be any number. Then for random $a,b$ the probability that $d \mid a+b$ is $1/d$. That means that if $v,w \in \mathbb{Z}^{1\times n}$ are vectors, then $d \mid v+w$ with a probability of $1/d^n$. This effect is noticeable in particular for small numbers like $d=2,3$ and in the last iterations of the $LD^{-1}U$ decomposition when the number of non-zero entries in the rows has shrunk. For instance, in the second last iterations we only have three rows with at most three non-zero entries each. Moreover, we know that the first non-zero entries of the rows cancel during cross-multiplication. Thus, a factor of 2 appears with a probability of 25% in one of those rows, a factor of 3 with a probability of 11.11%. In the example above, the probability for the factor 5 to appear in the fourth row was 4%.

## 5. Expected Number of Factors

In this section, we provide a detailed analysis of the expected number of common "statistical" factors in the rows of $U$, in the case when the input matrix $A$ has integer entries, that is, $\mathbb{D}=\mathbb{Z}$. We base our considerations on a "uniform" distribution on $\mathbb{Z}$, e.g., by imposing a uniform distribution on $\{-n,\ldots,n\}$ for very large $n$. However, the only relevant property that we use is the assumption that the probability that a randomly chosen integer is divisible by $p$ is $1/p$.

We consider a matrix $A = (A_{i,j})_{1\le i,j\le n} \in \mathbb{Z}^{n\times n}$ of full rank. The assumption that $A$ be square is made for the sake of simplicity; the results shown below immediately generalize to rectangular matrices. As before, let $U$ be the upper triangular matrix from the $LD^{-1}U$ decomposition of $A$:

$$U = \begin{pmatrix} U_{1,1} & U_{1,2} & \cdots & U_{1,n} \\ 0 & U_{2,2} & \cdots & U_{2,n} \\ \vdots & & \ddots & \vdots \\ 0 & \cdots & & U_{n,n} \end{pmatrix}.$$

Define

$$g_k := \gcd(U_{k,k}, U_{k,k+1}, \ldots, U_{k,n})$$

to be the greatest common divisor of all entries in the $k^{\text{th}}$ row of $U$. Counting (with multiplicities) all the prime factors of $g_1,\ldots,g_{n-1}$, one gets the picture shown in Figure 1; $g_n$ is omitted as it contains only the single nonzero entry $U_{n,n}=\det(A)$. Our goal is to give a probabilistic explanation for the occurrence of these common factors, whose number seems to grow linearly with the dimension of the matrix.

As we have seen in the proof of Theorem 6, the entries $U_{k,\ell}$ can be expressed as minors of the original matrix $A$:

$$U_{k,\ell} = \det \begin{pmatrix} A_{1,1} & A_{1,2} & \cdots & A_{1,k-1} & A_{1,\ell} \\ A_{2,1} & A_{2,2} & \cdots & A_{2,k-1} & A_{2,\ell} \\ \vdots & \vdots & & \vdots & \vdots \\ A_{k,1} & A_{k,2} & \cdots & A_{k,k-1} & A_{k,\ell} \end{pmatrix}.$$

Observe that the entries $U_{k,\ell}$ in the $k^{\text{th}}$ row of $U$ are all given as determinants of the same matrix, where only the last column varies. For any integer $q \ge 2$ we have that $q \mid g_k$ if $q$ divides all these determinants. A sufficient condition for the latter to happen is that the determinant

$$h_k := \det \begin{pmatrix} A_{1,1} & \dots & A_{1,k-1} & 1 \\ A_{2,1} & \dots & A_{2,k-1} & x \\ \vdots & \vdots & & \vdots \\ A_{k,1} & \dots & A_{k,k-1} & x^{k-1} \end{pmatrix}$$

is divisible by $q$ as a polynomial in $\mathbb{Z}[x]$, i.e., if $q$ divides the content of the polynomial $h_k$. We now aim at computing how likely it is that $q \mid h_k$ when $q$ is fixed and when the matrix entries $A_{1,1}, \dots, A_{k,k-1}$ are chosen randomly. Since $q$ is now fixed, we can equivalently study this problem over the finite ring $\mathbb{Z}_q$, which means that the matrix entries are picked randomly and uniformly from the finite set $\{0, \dots, q-1\}$. Moreover, it turns out that it suffices to answer this question for prime powers $q = p^j$.

The probability that all $k \times k$-minors of a randomly chosen $k \times (k+1)$-matrix are divisible by $p^j$, where $p$ is a prime number and $j \geq 1$ is an integer, is given by

$$P_{p,j,k} := 1 - \left(1 + p^{1-j-k} \frac{p^k - 1}{p - 1}\right) \prod_{i=0}^{k-1} (1 - p^{-j-i}),$$

which is a special case of Brent and McKay [2, Thm. 2.1]. Note that this is exactly the probability that $h_{k+1}$ is divisible by $p^j$. Recalling the definition of the $q$-Pochhammer symbol

$$(a; q)_k := \prod_{i=0}^{k-1} (1 - aq^i), \quad (a; q)_0 := 1,$$

the above formula can be written more succinctly as

$$P_{p,j,k} := 1 - \left(1 + p^{1-j-k} \frac{p^k - 1}{p - 1}\right) \left(\frac{1}{p^j}; \frac{1}{p}\right)_k.$$

Now, an interesting observation is that this probability does not, as one could expect, tend to zero as $k$ goes to infinity. Instead, it approaches a nonzero constant that depends on $p$ and $j$ (see Table 1):

$$P_{p,j,\infty} := \lim_{k \to \infty} P_{p,j,k} = 1 - \left(1 + \frac{p^{1-j}}{p-1}\right) \left(\frac{1}{p^j}; \frac{1}{p}\right)_\infty$$

| $p^j$ | $k=1$ | $k=2$ | $k=3$ | $k=4$ | $k=5$ | $k=6$ | $k=\infty$ |
|---|---|---|---|---|---|---|---|
| 2 | 0.25000 | 0.34375 | 0.38477 | 0.40399 | 0.41330 | 0.41789 | 0.42242 |
| 3 | 0.11111 | 0.14403 | 0.15460 | 0.15808 | 0.15923 | 0.15962 | 0.15981 |
| 4 | 0.06250 | 0.09766 | 0.11560 | 0.12461 | 0.12912 | 0.13138 | 0.13364 |
| 5 | 0.04000 | 0.04768 | 0.04920 | 0.04951 | 0.04957 | 0.04958 | 0.04958 |
| 7 | 0.02041 | 0.02326 | 0.02367 | 0.02373 | 0.02374 | 0.02374 | 0.02374 |
| 8 | 0.01563 | 0.02588 | 0.03149 | 0.03440 | 0.03588 | 0.03662 | 0.03737 |

TABLE 1. Behavior of the sequence $(P_{p,j,k})_{k \in \mathbb{N}}$ for some small values of $p^j$.

Using the probability $P_{p,j,k}$, one can write down the expected number of factors in the determinant $h_{k+1}$, i.e., the number of prime factors in the content of the polynomial $h_{k+1}$, counted with multiplicities:

$$\sum_{p \in \mathbb{P}} \sum_{j=1}^{\infty} P_{p,j,k},$$

where $\mathbb{P} = \{2, 3, 5, \dots\}$ denotes the set of prime numbers. The inner sum can be simplified as follows, yielding the expected multiplicity $M_{p,k}$ of a prime factor $p$ in $h_{k+1}$:

$$M_{p,k} := \sum_{j=1}^{\infty} P_{p,j,k} = \sum_{j=1}^{\infty} \left(1 - \left(1 + p^{1-j-k}\frac{p^k-1}{p-1}\right)\left(\frac{1}{p^j}; \frac{1}{p}\right)_k\right)$$

$$= -\sum_{j=1}^{\infty}\left(\left(\frac{1}{p^j};\frac{1}{p}\right)_k - 1\right) - p^{1-k}\frac{p^k-1}{p-1}\sum_{j=1}^{\infty}\frac{1}{p^j}\left(\frac{1}{p^j};\frac{1}{p}\right)_k$$

$$= -\sum_{j=1}^{\infty}\sum_{i=1}^{k}(-1)^i p^{-ij-i(i-1)/2}\begin{bmatrix}k\\i\end{bmatrix}_{1/p} - p^{1-k}\frac{p^k-1}{p-1}\frac{p^k}{p^{k+1}-1}$$

$$= \sum_{i=1}^{k}\frac{(-1)^{i-1}}{p^{i(i-1)/2}(p^i-1)}\begin{bmatrix}k\\i\end{bmatrix}_{1/p} + \frac{1}{p^{k+1}-1} - \frac{1}{p-1}$$

In this derivation we have used the expansion formula of the $q$-Pochhammer symbol in terms of the $q$-binomial coefficient

$$\begin{bmatrix}n\\k\end{bmatrix}_q := \frac{(1-q^n)(1-q^{n-1})\cdots(1-q^{n-k+1})}{(1-q^k)(1-q^{k-1})\cdots(1-q)},$$

evaluated at $q = 1/p$. Moreover, the identity that is used in the third step,

$$\sum_{j=1}^{\infty}\frac{1}{p^j}\left(\frac{1}{p^j};\frac{1}{p}\right)_k = \frac{p^k}{p^{k+1}-1},$$

is certified by rewriting the summand as

$$\frac{1}{p^j}\left(\frac{1}{p^j};\frac{1}{p}\right)_k = t_{j+1} - t_j \quad \text{with} \quad t_j = \frac{p^k(p^{1-j}-1)}{p^{k+1}-1}\left(\frac{1}{p^j};\frac{1}{p}\right)_k$$

and by applying a telescoping argument.

Hence, when we let $k$ go to infinity, we obtain

$$M_{p,\infty} = \lim_{k\to\infty}\sum_{j=1}^{\infty} P_{p,j,k} = \sum_{i=1}^{\infty}\frac{(-1)^{i-1}}{p^{i(i-1)/2}(p^i-1)}\frac{(p^{-i-1};p^{-1})_\infty}{(p^{-1};p^{-1})_\infty} - \frac{1}{p-1}.$$

Note that the sum converges quickly, so that one can use the above formula to compute an approximation for the expected number of factors in $h_{k+1}$ when $k$ tends to infinity

$$\sum_{p\in\mathbb{P}} M_{p,\infty} \approx 0.89764,$$

which gives the asymptotic slope of the function plotted in Figure 1.

As discussed before, the divisibility of $h_k$ by some number $q \geq 2$ implies that the greatest common divisor $g_k$ of the $k^{\text{th}}$ row is divisible by $q$, but this is not a necessary condition. It may happen that $h_k$ is not divisible by $q$, but nevertheless $q$ divides each $U_{k,\ell}$ for $k \leq \ell \leq n$. The probability for this to happen is the same as the probability that the greatest common divisor of $n-k+1$ randomly chosen integers is divisible by $q$. The latter obviously is $q^{-(n-k+1)}$. Thus, in addition to the factors coming from $h_k$, one can expect

$$\sum_{p\in\mathbb{P}}\sum_{j=1}^{\infty}\frac{1}{p^{j(n-k+1)}} = \sum_{p\in\mathbb{P}}\frac{1}{p^{n-k+1}-1}$$

many prime factors in $g_k$.

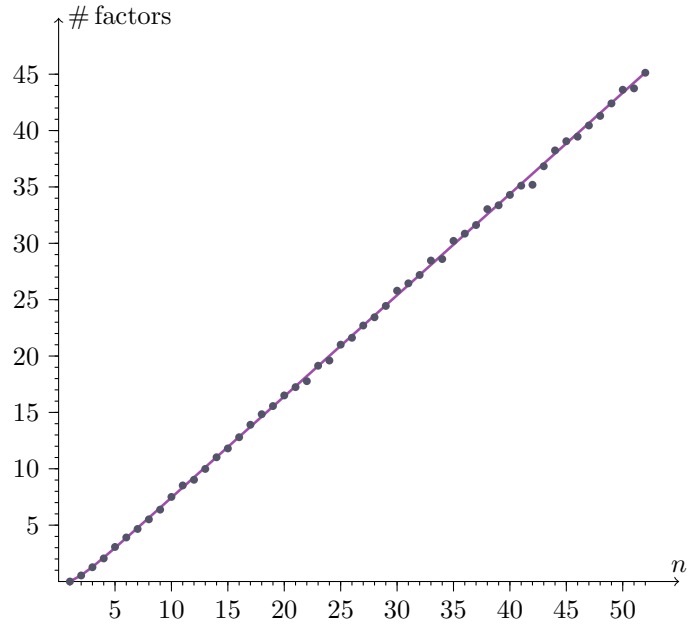

FIGURE 1. Number of factors depending on the size $n$ of the matrix. The curve shows the function $F(n)$, while the dots represent experimental data: for each dimension $n$, 1000 matrices were generated with random integer entries between 0 and $10^9$.

Summarizing, the expected number of prime factors in the rows of the matrix $U$ is

$$F(n) = \sum_{k=2}^{n-1} \sum_{p \in \mathbb{P}} M_{p,k-1} + \sum_{k=1}^{n-1} \sum_{p \in \mathbb{P}} \frac{1}{p^{n-k+1} - 1}$$

$$= \sum_{p \in \mathbb{P}} \left( \sum_{k=0}^{n-2} M_{p,k} + \sum_{k=0}^{n-2} \frac{1}{p^{k+2} - 1} \right)$$

$$= \sum_{p \in \mathbb{P}} \sum_{k=0}^{n-2} \left( \sum_{i=1}^{k} \frac{(-1)^{i-1}}{p^{i(i-1)/2}(p^i - 1)} \begin{bmatrix} k \\ i \end{bmatrix}_{1/p} + \frac{1}{p^{k+2} - 1} + \frac{1}{p^{k+1} - 1} - \frac{1}{p - 1} \right).$$

¿From the discussion above, it follows that for large $n$ this expected number can be approximated by a linear function as follows:

$$F(n) \approx 0.89764\, n - 1.53206.$$

## 6. $QR$ Decomposition

A fraction-free $QR$ decomposition, which is based on the $LD^{-1}U$ decomposition, was given in Zhou and Jeffrey [18]. In this section, we present a refined version of this algorithm (see Theorem 12). As a first step in its proof, we will need the Cholesky decomposition, which is introduced in the following lemma.

This section assumes that $\mathbb{D}$ has characteristic 0 which is needed in order to assure that $A^t A$ has full rank.

**Theorem 11.** *Let $A \in \mathbb{D}^{n \times n}$ be a symmetric matrix such that its $LD^{-1}U$ decomposition can be computed without permutations; then we have $U = L^t$, that is,*

$$A = LD^{-1}L^t.$$

*Proof.* Compute the decomposition $A = LD^{-1}U$ as in Theorem 1. If we do not execute item 4 of Algorithm 2, we obtain the decomposition

$$A = \tilde{L}\tilde{D}^{-1}\tilde{U} = \begin{pmatrix} \mathcal{L} & \mathbf{0} \\ \mathcal{M} & \mathbf{1} \end{pmatrix} \begin{pmatrix} D & \mathbf{0} \\ \mathbf{0} & \mathbf{1} \end{pmatrix}^{-1} \begin{pmatrix} \mathcal{U} & \mathcal{V} \\ \mathbf{0} & \mathbf{0} \end{pmatrix}.$$

Then because $A$ is symmetric, we obtain

$$\tilde{L}\tilde{D}^{-1}\tilde{U} = A = A^t = \tilde{U}^t\tilde{D}^{-1}\tilde{L}^t$$

The matrices $\tilde{L}$ and $\tilde{D}$ have full rank which implies

$$\tilde{U}(\tilde{L}^t)^{-1}\tilde{D} = \tilde{D}\tilde{L}^{-1}\tilde{U}^t.$$

Examination of the matrices on the left hand side reveals that they are all upper triangular. Therefore also their product is an upper triangular matrix. Similarly, the right hand side is a lower triangular matrix and the equality of the two implies that they must both be diagonal. Canceling $\tilde{D}$ and rearranging the equation yields $\tilde{U} = (\tilde{L}^{-1}\tilde{U}^t)\tilde{L}^t$ where $\tilde{L}^{-1}\tilde{U}^t$ is diagonal. This shows that the rows of $\tilde{U}$ are just multiples of the rows of $\tilde{L}^t$. However, we know that the first $r$ diagonal entries of $\tilde{U}$ and $\tilde{L}$ are the same, where $r$ is the rank of $\tilde{U}$. This yields

$$\tilde{L}^{-1}\tilde{U}^t = \begin{pmatrix} \mathbf{1}_r & \mathbf{0} \\ \mathbf{0} & \mathbf{0} \end{pmatrix},$$

and hence, when we remove the unnecessary last $n - r$ rows of $\tilde{U}$ and the last $n - r$ columns of $\tilde{L}$ (as suggested in Jeffrey [10]), we remain with $U = L^t$. □

The following theorem is a variant of Zhou and Jeffrey [18, Thm. 8], where we exploit the symmetry of $A^t A$ by invoking Theorem 11. This leads to a nicer representation of the decomposition, and we obtain more information about $\Theta^t\Theta$.

**Theorem 12.** *Let $A \in \mathbb{D}^{m \times n}$ with $n \leq m$ and with full column rank. Then the partitioned matrix $(A^t A \mid A^t)$ has $LD^{-1}U$ decomposition*

$$(A^t A \mid A^t) = R^t D^{-1}(R \mid \Theta^t),$$

*where $\Theta^t\Theta = D$ and $A = \Theta D^{-1}R$.*

*Proof.* Since $A$ has full column rank, the *Gramian matrix* $A^t A$ will have full rank, too. By taking the first $k$ columns of $A$ (and the first $k$ rows of $A^t$), it follows that also the $k^{\text{th}}$ principal minor of $A^t A$ is nonzero. Consequently, when we compute the $LD^{-1}U$ decomposition, we do not need any permutations.

Hence, by Theorem 11, we can decompose the symmetric matrix $A^t A$ as

$$A^t A = R^t D^{-1} R.$$

Applying the same row transformations to $A^t$ yields a matrix $\Theta^t$, that is, we obtain $(A^t A \mid A^t) = R^t D^{-1}(R \mid \Theta^t)$. As in the proof of Zhou and Jeffrey [18, Thm. 8], we easily compute that $A = \Theta D^{-1}R$ and that $\Theta^t\Theta = D^t(R^{-1})^t A^t A R^{-1} D = D^t(R^{-1})^t R^t D^{-1} R R^{-1} D = D$. □

For example, let $A \in \mathbb{Z}[x]^{3 \times 3}$ be the matrix

$$A = \begin{pmatrix} x & 1 & 2 \\ 2 & 0 & -x \\ x & 1 & x+1 \end{pmatrix}.$$

Then the $LD^{-1}U$ decomposition of $A^tA = R^tD^{-1}R$ is given by

$$R = \begin{pmatrix} 2(x^2+2) & 2x & x(x+1) \\ 0 & 8 & 4(x^2+x+3) \\ 0 & 0 & 4(x-1)^2 \end{pmatrix},$$

$$D = \begin{pmatrix} 2(x^2+2) & 0 & 0 \\ 0 & 16(x^2+2) & 0 \\ 0 & 0 & 32(x-1)^2 \end{pmatrix},$$

and we obtain for the $QR$ decomposition $A = \Theta D^{-1}R$:

$$\Theta = \begin{pmatrix} x & 4 & -4(x-1) \\ 2 & -4x & 0 \\ x & 4 & 4(x-1) \end{pmatrix}.$$

We see that the $\Theta D^{-1}R$ decomposition has some common factor in the last column of $\Theta$. This observation is explained by the following theorem.

**Theorem 13.** *With full-rank $A \in \mathbb{D}^{n \times n}$ and $\Theta$ as in Theorem 12, we have for all $i = 1, \ldots, n$ that*

$$\Theta_{in} = (-1)^{n+i} \det_{i,n} A \cdot \det A$$

*where $\det_{i,n} A$ is the $(i,n)$ minor of $A$.*

*Proof.* We use the notation from the proof of Theorem 12. From $\Theta D^{-1}R = A$ and $\Theta^t\Theta = D$ we obtain

$$\Theta^t A = \Theta^t \Theta D^{-1} R = R.$$

Thus, since $A$ has full rank, $\Theta^t = RA^{-1}$ or, equivalently,

$$\Theta = (RA^{-1})^t = (A^{-1})^t R^t = (\det A)^{-1}(\operatorname{adj} A)^t R^t$$

where $\operatorname{adj} A$ is the adjoint matrix of $A$. Since $R^t$ is a lower triangular matrix with $\det A^t A = (\det A)^2$ at position $(n,n)$, the claim follows. □

For the other columns of $\Theta$ we can state the following.

**Theorem 14.** *The $k^{th}$ determinantal divisor $d_k^*$ of $A$ divides the $k^{th}$ column of $\Theta$ and the $k^{th}$ row of $R$. Moreover, $d_{k-1}^* d_k^*$ divides $D_{k,k}$ for $k \geq 2$.*

*Proof.* We first show that the $k^{th}$ determinantal divisor $\delta_k^*$ of $(A^tA \mid A^t)$ is the same as $d_k^*$. Obviously, $\delta_k^* \mid d_k^*$ since all minors of $A$ are also minors of the right block $A^t$ of $(A^tA \mid A^t)$. Consider now the left block $A^tA$. We have by the Cauchy–Binet theorem [3, § 4.6]

$$\det_{I,J}(A^tA) = \sum_{\substack{K \subseteq \{1, \ldots, n\} \\ |K| = q}} (\det_{K,I} A)(\det_{K,J} A)$$

where $I, J \subseteq \{1, \ldots, n\}$ with $|I| = |J| = q \geq 1$ are two index sets and $\det_{I,J} M$ denotes the minor for these index sets of a matrix $M$. Thus, $(d_k^*)^2$ divides any minor of $A^tA$ since it divides every summand on the right hand side; and we see that $d_k^* \mid \delta_k^*$.

Now, we use Theorem 12 and Theorem 6 to conclude that $d_k^*$ divides the $k^{th}$ row of $(R \mid \Theta^t)$ and hence the $k^{th}$ row of $R$ and the $k^{th}$ column of $\Theta$. Moreover, $D_{k,k} = R_{k-1,k-1}R_{k,k}$ for $k \geq 2$ by Theorem 1 which implies $d_{k-1}^* d_k^* \mid D_{k,k}$. □

Knowing that there is always a common factor, we can cancel it, which leads to a fraction-free $QR$ decomposition of smaller size.

**Theorem 15.** *Given a square matrix $A$, a reduced fraction-free $QR$ decomposition is given by $A = \hat{\Theta}\hat{D}^{-1}\hat{R}$, where $S = \operatorname{diag}(1, 1, \ldots, \det A)$ and $\hat{\Theta} = \Theta S^{-1}$, and $\hat{R} = S^{-1}R$. In addition, $\hat{D} = S^{-1}DS^{-1} = \hat{\Theta}^t\hat{\Theta}$.*

*Proof.* By Theorem 13, $\Theta S^{-1}$ is an exact division. The statement of the theorem follows from $A = \Theta S^{-1} S D^{-1} S S^{-1} R$. □

If we apply Theorem 15 to our previous example, we obtain the simpler $QR$ decomposition, where the factor $\det A = -2(x-1)$ has been removed.

$$\begin{pmatrix} x & 4 & 2 \\ 2 & -4x & 0 \\ x & 4 & -2 \end{pmatrix} \begin{pmatrix} 2(x^2+2) & 0 & 0 \\ 0 & 16(x^2+2) & 0 \\ 0 & 0 & 8 \end{pmatrix}^{-1} \begin{pmatrix} 2(x^2+2) & 2x & x(x+1) \\ 0 & 8 & 4(x^2+x+3) \\ 0 & 0 & -2(x-1) \end{pmatrix}.$$

The properties of the $QR$-decomposition are strong enough to guarantee a certain uniqueness of the output.

**Theorem 16.** *Let $A \in \mathbb{D}^{n \times n}$ have full rank. Let $A = \Theta D^{-1} R$ the decomposition from Theorem 12; and let $A = \tilde{\Theta} \tilde{D}^{-1} \tilde{R}$ be another decomposition where $\tilde{\Theta}, \tilde{D}, \tilde{R} \in \mathbb{D}^{n \times n}$ are such that $\tilde{D}$ is a diagonal matrix, $\tilde{R}$ is an upper triangular matrix and $\tilde{\Theta}^t \tilde{\Theta}$ is a diagonal matrix. Then $\Theta^t \tilde{\Theta}$ is also a diagonal matrix and $\tilde{R} = (\Theta^t \tilde{\Theta})^{-1} \tilde{D} R$.*

*Proof.* We have
$$\tilde{\Theta} \tilde{D}^{-1} \tilde{R} = \Theta D^{-1} R \qquad \text{and thus} \qquad \Theta^t \tilde{\Theta} \tilde{D}^{-1} \tilde{R} = \Theta^t \Theta D^{-1} R = R.$$
If $R$ and $\tilde{R}$ have full rank, this is equivalent to
$$\Theta^t \tilde{\Theta} = R \tilde{R}^{-1} \tilde{D}.$$
Note that all the matrices on the right hand side are upper triangular. Similarly, we can compute that
$$\tilde{\Theta}^t \Theta D^{-1} R = \tilde{\Theta}^t \tilde{\Theta} \tilde{D}^{-1} \tilde{R} = \Delta \tilde{D}^{-1} \tilde{R}$$
which implies $\tilde{\Theta}^t \Theta = \Delta \tilde{D}^{-1} \tilde{R} R^{-1} D$. Hence, also $\tilde{\Theta}^t \Theta = (\Theta^t \tilde{\Theta})^t$ is upper triangular and consequently $\tilde{\Theta}^t \Theta = T$ for some diagonal matrix $T$ with entries from $\mathbb{D}$. We obtain $R = T \tilde{D}^{-1} \tilde{R}$ and thus $\tilde{R} = T^{-1} \tilde{D} R$. □

## 7. Acknowledgments

We would like to thank Kevin G. Hare and Arne Winterhof for helpful comments and discussions. We are grateful to the anonymous referees whose insightful remarks improved this paper considerably. In particular, the statements of Corollary 8 and Theorem 14 were pointed out by one of the referees.

## References

[1] Bareiss, E. H.: 1968, Mathematics of Computation 22(103), 565
[2] Brent, R. P. and McKay, B. D.: 1987, Discrete Mathematics pp 35–49
[3] Broida, J. G. and Williamson, S. G.: 1989, A Comprehensive Introduction to Linear Algebra, Addison Wesley
[4] Corless, R. M. and Jeffrey, D. J.: 1997, SIGSAM Bulletin 31(3), 20
[5] Erlingsson, Ú., Kaltofen, E., and Musser, D.: 1996, in International Symposium on Symbolic and Algebraic Computation, pp 275–282, ACM Press
[6] Geddes, K. O., Czapor, S. R., and Labahn, G.: 1992, Algorithms for Computer Algebra, Kluwer Academic Publisher
[7] Giesbrecht, M. W. and Storjohann, A.: 2002, Journal of Symbolic Computation 34(3), 157
[8] Hare, K. G.: 2016, Personal Communication
[9] Jebelean, T.: 1993, J. Symbolic Computation 15, 169
[10] Jeffrey, D. J.: 2010, ACM Communications in Computer Algebra 44(171), 1
[11] Kaltofen, E. and Yuhasz, G.: 2013, Linear Algebra and its Applications 439(9), 2515
[12] Lee, H. R. and Saunders, B. D.: 1995, Journal of Symbolic Computation 19(5), 393

[13] Middeke, J. and Jeffrey, D. J.: 2014, Fraction-free factoring revisited, Poster presentation at the International Symposium on Symbolic and Algebraic Computation

[14] Nakos, G. C., Turner, P. R., and Williams, R. M.: 1997, SIGSAM Bulletin 31(3), 11

[15] Newman, M.: 1972, Integral Matrices, Vol. 45 of Pure and Applied Mathematics, Academic Press, New York

[16] Pauderis, C. and Storjohann, A.: 2013, in M. Kauers (ed.), Proceedings of the International Symposium on Symbolic and Algebraic Computation, pp 307–314, ACM Press

[17] Winterhof, A.: 2016, Personal Communication

[18] Zhou, W. and Jeffrey, D. J.: 2008, Frontiers of Computer Science in China 2(1), 67

Johannes Middeke
Research Institute for Symbolic Computation, Johannes Kepler University, Altenberger Straße 69, A-4040 Linz, Austria
e-mail: jmiddeke@risc.jku.at

David J. Jeffrey
Department of Applied Mathematics, University of Western Ontario, Middlesex College, Room 255, 1151 Richmond Street North, London, Ontario, Canada, N6A 5B7
e-mail: djeffrey@uwo.ca

Christoph Koutschan
Johann Radon Institute for Computational and Applied Mathematics, Austrian Academy of Sciences, Altenberger Straße 69, A-4040 Linz, Austria
e-mail: christoph.koutschan@ricam.oeaw.ac.at